\documentclass[reqno]{amsart}
\usepackage{basicmacros}

\title[Entropy for Non-Uniformly Continuous Maps]{Topological Entropy\\ for Non-Uniformly Continuous Maps
}
\author{Boris Hasselblatt}
	\address{Department of Mathematics, Tufts University, Medford, MA 02155}
	\email{boris.hasselblatt@tufts.edu}
\author{Zbigniew Nitecki}
	\address{Department of Mathematics, Tufts University, Medford, MA 02155}
	\email{zbigniew.nitecki@tufts.edu}
\author{James Propp}
	\address{Department of Mathematics, University of Wisconsin, Madison, WI 53706}
	\email{propp@math.wisc.edu}
\date{\today} 

\newcommand{\dist}{\ensuremath{d}}
\newcommand{\dis}[2]{\ensuremath{\dist(#1,#2)}}
\newcommand{\dists}[1]{\ensuremath{\dist_{#1}}}
\newcommand{\diss}[3]{\ensuremath{\dists{#1}(#2,#3)}}
\newcommand{\dBowen}[2]{\ensuremath{\dist^{#1}_{#2}}}
\newcommand{\dB}[4]{\ensuremath{\dBowen{#1}{#2}(#3,#4)}}
\newcommand{\dnf}{\dBowen{f}{n}}
\newcommand{\dnfof}[2]{\dnf({#1},{#2})}
\newcommand{\tdist}{\ensuremath{\tilde{d}}}
\newcommand{\tdis}[2]{\ensuremath{\tdist({#1},{#2})}}
\newcommand{\tdBowen}[2]{\ensuremath{\tdist^{#1}_{#2}}}
\newcommand{\tdB}[4]{\ensuremath{\tdBowen{#1}{#2}(#3,#4)}}

\newcommand{\ent}{\ensuremath{\mathfrak{h}}}
\newcommand{\hsub}[1]{\ensuremath{\ent_{#1}}}
\newcommand{\htop}[1]{\ensuremath{\ent_{top}(#1)}}

\newcommand{\hsep}{\hsub{BD}}
\newcommand{\hsepof}[3]{\ensuremath{\hsep(#1,#2,#3)}}
\newcommand{\hcomp}{\hsub{Bc}}
\newcommand{\hcompof}[3]{\ensuremath{\hcomp(#1,#2,#3)}}

\newcommand{\diam}{\operatorname{diam}}

\renewcommand{\ms}[3]{\ensuremath{{\operatorname{maxsep}}[#1,#2,#3]}}
\newcommand{\minsp}[3]{\ensuremath{{\operatorname{minspan}}[#1,#2,#3]}}
\newcommand{\minspan}[3]{\ensuremath{{\operatorname{minspan}}[#1,#2,#3]}}
\newcommand{\GR}[1]{\ensuremath{GR\{#1\}}}
\newcommand{\sep}[2]{$(#1,#2)$-separated}

\newcommand{\nepssep}{\sep{n}{\eps}}

\newcommand{\seqx}{\ensuremath{\mathbf{x}}}
\newcommand{\seqxp}{\ensuremath{\seqx^\prime}}

\newcommand*{\Xs}[1]{\ensuremath{\X_{#1}}}
\newcommand{\tilf}{\ensuremath{\tilde{f}}}
\newcommand{\tX}{\ensuremath{\tilde{\X}}} 
\newcommand{\tfXX}{\selfmap{\tilf}{\tX}}
\newcommand{\Xd}{\pairord{\X}{\dist}}
\newcommand{\Xdnf}{\pairord{\X}{\dnf}}
\newcommand{\fXd}{\selfmap{\f}{\Xd}}
\newcommand{\tXd}{\pairord{\tX}{\tdist}}
\newcommand{\tfXd}{\selfmap{\tilf}{\tXd}}

\newcommand{\CP}[1]{\ensuremath{\mathbb{CP}^{#1}}}

\renewcommand{\ent}{\ensuremath{h}}% change from frak in entropymacros.tex

\newcommand{\hX}{\ensuremath{\hat{\X}}}
\newcommand{\hdist}{\ensuremath{\hat{\dist}}}
\newcommand{\hdis}[2]{\ensuremath{\hdist(#1,#2)}}
\newcommand{\hXd}{\pairord{\hX}{\hdist}}

\newcommand{\hXfd}{\pairord{\Xf}{\hdist}}

\newcommand{\XtoN}{\ensuremath{\X^{\Nat}}}
\newcommand{\Xf}{\ensuremath{\X^{f}}}

\newcommand{\cdist}{\hdist}

\newcommand{\hF}{\hsub{F}}
\newcommand{\hFof}[1]{\ensuremath{\hF(#1)}}
\newcommand{\hFf}{\hFof{\f}}

\newcommand{\hAKMof}[1]{\ensuremath{\hsub{AKM}(#1)}}
\newcommand{\hAKMf}{\hAKMof{f}}

\newcommand{\Hof}[1]{\ensuremath{H\left(#1\right)}}

\newcommand{\peps}{\ensuremath{\eps\pr}}
\newcommand{\spr}{\ensuremath{s\pr}}

\newcommand{\discD}{\ensuremath{\mathbb{D}}}
\newcommand{\disc}{\ensuremath{\discD}}
\newcommand{\pdisc}{\ensuremath{\discD^{\circ}}}

\newcommand{\Selfmap}{self-map}

\newcommand{\dn}{\ensuremath{d_{N}}}

\newcommand{\finv}{\ensuremath{\inverse{f}}}

\begin{document}
\begin{abstract}
	The literature contains several extensions of the standard definitions of topological entropy 
	for a continuous
	\Selfmap{} \fXX{} from the case when \X{} is a compact metric space to the case when \X{} is 
	allowed to be non-compact.  These extensions all require the space \X{} to be totally bounded,
	or equivalently to have a compact completion, and are invariants of uniform conjugacy.
	When the map \f{} is uniformly continuous, it 
	extends continuously to the completion, 
	and the various notions of entropy reduce to the standard ones
	(applied to this extension).  However, when uniform continuity is not assumed, 
	these new quantities
	can differ.  We consider extensions proposed by Bowen 
	(maximizing over compact subsets)
	and Friedland (using the compactification of the graph of \f{}) 
	as well as a straightforward extension
	of Bowen and Dinaburg's definition from the compact case, 
	assuming that \X{} is totally bounded, 
	but not necessarily compact.  
	This
	last extension agrees with Friedland's, and both dominate the one proposed by Bowen
	(Theorem 6).  
	Examples show how varying the metric outside its uniform class can vary both
	quantities.  The natural
	extension of Adler-Konheim-McAndrews' original (metric-free) definition of topological entropy
	beyond compact spaces is unfortunately infinite for a great number of noncompact examples
	(Proposition 7). 
	
\end{abstract} 

\maketitle

There are two standard definitions of topological entropy for a continuous \Selfmap{} of a compact
metric space.  The original definition by Adler, Konheim and McAndrew \cite{Adler:top}, based on open covers, can in principle be applied to a continuous \Selfmap{}
of any compact topological space, while the reformulation of this definition by Bowen \cite{Bowen:homog} and Dinaburg \cite{Dinaburg}, based on 
the dispersion of orbits, requires a metric. When the metric space is compact,
these two definitions yield the same quantity, which is an invariant of topological conjugacy.
In particular,  the Bowen-Dinaburg version of entropy is independent of the (compact) metric used to compute it.

When \X{} is not compact, the situation is more complicated, as a number of invariants that are always equal in the compact case can differ in a noncompact setting.  To obtain a nontrivial invariant, one must take steps to preserve some features associated with compactness.  In \cite{Bowen:homog}, Bowen proposed an invariant based on measuring the dispersion of orbits emanating from a compact subset $K\subset\X$ and taking the supremum over all such subsets $K$.   Bowen was motivated by uniformly continuous examples, and his definition has been taken as a ``standard'' definition
of topological entropy for a (uniformly continuous) \Selfmap{} of a metric space which is not assumed to be compact (see, for example, \cite[pp. 168-176]{Walters}). We shall formulate this in \refer{sec}{noncompact} and call it \emph{Bowen compacta entropy}.  

However, in many cases of interest, such as billiards on tables with corners or meromorphic \Selfmap{}s on complex projective spaces, \X{} occurs naturally as a subset of a compact metrizable ambient space, giving a preferred class of metrics, but the map is not uniformly continuous in this metric. 
Friedland \cite{Friedland:polyrat}, motivated by examples of the second type, started from an interpretation of the Bowen-Dinaburg calculation (in the compact case) by Gromov \cite{Gromov:entholo}, and proposed a different invariant, based on a compactification of the graph of the map (Subsection 2.2).  

Both of these invariants begin from a metric which, as we shall see, must have a compact completion; 
both are unchanged if the metric is replaced by a uniformly equivalent one (so that the completions are homeomorphic) but both can change if the new metric is equivalent, but not uniformly equivalent, to the original one.

In this note we approach this situation abstractly and intrinsically.  We start with a continuous \Selfmap{} \fXX{}, assuming that \X{} is a (not necessarily closed) subset of a compact metric space. Thus, the restriction of the ambient metric to \X{} is totally bounded---but the map \f{} is not assumed to be uniformly continuous, and hence need not extend continuously to the (compact) closure of \X{} in the ambient space.   Our main observation is that the Bowen-Dinaburg calculation can be used \verbat{} in this context\footnote{Such a procedure is followed, without comment, in for example \cite{DinhSibony}.},
and yields the invariant defined by Friedland, which in turn dominates the invariant proposed by Bowen
(Theorem 6).
In Section 4, we consider some examples that illustrate the way the choice of metric can affect the values of both invariants\footnote{Walters gives an example \cite[p. 171]{Walters} of two equivalent, but not uniformly equivalent, metrics on the line for which his version of topological entropy---which is the same as the Bowen compacta entropy---for the map $x\mapsto2x$ has different values.}, 
and which also indicate the failure of a certain natural strategy for extending the original Adler-Konheim-McAndrew definition to this situation (Proposition 7).

\section{Entropy in Compact Spaces}\label{sec:compact}

We briefly review the definitions of Adler-Konheim-McAndrew and Bowen-Dinaburg, for future reference.

\subsection{Adler-Konheim-McAndrew's Definition}
In \cite{Adler:top}, the topological entropy of a continuous \Selfmap{} \fXX{} of a compact topological space is defined as follows.  Given an open cover \alf{} of \X{}, denote by \Hof{\alf} the logarithm of the cardinality of a minimal subcover.  The entropy of \f{} relative to \alf{} is
\begin{equation}\label{eqn:relent}
	\hof{f,\alf}\eqdef\lim_{n\to\infty}\frac{1}{n}\Hof{\alfs{n}}
\end{equation}
where \alfs{n} is the mutual refinement of the covers formed by taking \alf{} together with all of its preimages under the first $n-1$ iterates of \f{}.  The topological entropy of \f{} is defined as the supremum of the entropy relative to all open covers of \X{}, or equivalently (in this case) all \emph{finite} open covers of \X{}.  To distinguish this calculation from others, we refer to this quantity as the 
\deffont{Adler-Konheim-McAndrew entropy}, and denote it by \hAKMf{}:
\begin{equation}\label{eqn:AKM}
	\hAKMf\eqdef\sup_{\alf\text{ finite open cover of \X}}\hof{f,\alf}.
\end{equation}

\subsection{Bowen-Dinaburg's Definition}
The Bowen-Dinaburg definition \cite{Bowen:homog, Dinaburg} of topological entropy for a continuous \Selfmap{} \fXX{} of a compact metric space can be formulated as follows.  First, some terminology.
 Given any metric space \pairord{\X}{\dist}, we call a subset $S\subset\X$ 
	\deffont{$\boldsymbol{\eps}$-separated} with respect to \dist{} for some \epsgo{}
	if distinct points of $S$ are spaced at least \eps{} apart: 
	\begin{equation*}\label{eqn:epsep}
		s\neq s\pr\in S\imply\dis{s}{s\pr}\geq\eps.
	\end{equation*}
	A set $S\subset\X$ \deffont{$\boldsymbol{\eps}$-spans} a subset $K\subset\X$ if every point of $K$ 
	is within distance \eps{}
	of some point of $S$:
	\begin{equation*}\label{eqn:epsspan}
		\forall x\in K\  \exists s\in S \text{ such that } \dis{x}{s}<\eps.
	\end{equation*}
	(When $S\subset K$, we can say $S$ is ``\eps-dense'' in $K$.)
	If there exist finite \eps-spanning sets for \X{} with respect to \dist{}, then since an 
	\eps-separated set
	which is maximal with respect to inclusion is also \eps-dense in \X{}, for any set $K\subset\X$ 
	the numbers
		\begin{align*}\label{eqn:maxsep}
			\minsp{K}{\dist}{\eps}
				&\eqdef\min\setbld{\card S}{S\subset\X\text{ \eps-spans $K$ with respect to \dist}}\\
			\ms{K}{\dist}{\eps}
				&\eqdef\max\setbld{\card S}{S\subset K\text{ is \eps-separated with respect to \dist}}
		\end{align*}
	are both finite and satisfy
		\begin{equation}\label{eqn:maxsepminspan}
			\minsp{K}{\dist}{\eps}\leq\ms{K}{\dist}{\eps}\leq\minsp{K}{\dist}{\frac{\eps}{2}}
		\end{equation}
	by an easy application of the triangle inequality.
	
	Now suppose \fXX{} is continuous with respect to a compact metric \dist{} on \X{}.  
	We construct the \deffont{Bowen-Dinaburg metrics} \dnf{}, \inflist{n}{1,2} via
		\begin{equation}\label{eqn:BowenDina}
			\dB{f}{n}{x}{\xp}\eqdef\max_{0\leq i<n}\dis{\ftoof{i}{x}}{\ftoof{i}{\xp}}.
		\end{equation}
	For each $n$, \dnf{} is another compact metric on \X{}, so the considerations above apply to each
	of the numbers \minsp{K}{\dnf}{\eps} and \ms{K}{\dnf}{\eps}, \inflist{n}{0,1,2}, \epsgo{}.
	We say that a set \deffont{$\boldsymbol{(n,\eps)}$-spans} $K\subset\X$ \resp{is 
	\deffont{$\boldsymbol{(n,\eps)}$-separated}} 
	if it \eps-spans $K$ \resp{is \eps-separated} with respect to \dnf{}.
	
	The (exponential) growth rate of any sequence \single{\cs{n}} of positive reals is defined by
		$$\GR{\cs{n}}\eqdef\limsup\frac{1}{n}\log\cs{n}.$$
		
	Then the entropy of \fXX{} on $K\subset\X$ with respect to the metric \dist{} is defined by taking
	the growth rate of the numbers \minsp{K}{\dnf}{\eps} or \ms{K}{\dnf}{\eps}
	 for \epsgo{} fixed (the two growth rates are related by \refer{eqn}{maxsepminspan}), 
	 then taking the limit as $\eps\to0$ (the two limits are equal by \refer{eqn}{maxsepminspan}). 
	  To distinguish this calculation for future reference, we will call it the 
	  \deffont{Bowen-Dinaburg entropy} of \f{} on $K$:
		\begin{align}\label{eqn:sepent}
			\hsepof{f}{K}{\dist}&\eqdef\lim_{\eps\to0}\GR{\ms{K}{\dnf}{\eps}}\\
				&=\lim_{\eps\to0}\GR{\minsp{K}{\dnf}{\eps}}.
		\end{align}
In particular, when \X{} is compact with respect to the metric \dist{}, then 
\begin{equation*}\label{eqn:sepAKM}
	\hsepof{f}{\X}{\dist}=\hAKMf.
\end{equation*}
(For more details, see \cite[Chap. 7]{Walters}.)

\subsection{Gromov's Observation}
Finally, we note an observation by Gromov \cite{Gromov:entholo} concerning the preceding definition.
(We formulate this in slightly different language than that used in \cite{Gromov:entholo}.)

The set of sequences in a compact space \X{} can be regarded as the product 
	$$\XtoN\eqdef\prod_{i=0}^{\infty}\Xs{i}=\setbld{\seqx\eqdef\xs{0}\xs{1}...}{\xs{i}\in\Xs{i}=\X}$$
of a sequence of copies of \X{} with the (Tikhonov) product topology.  When \X{} is a compact metric space, then so is \XtoN{}, and it is metrizable via 
the metric \hdist{} defined for any choice of $\rho>1$ by
\begin{equation*}\label{eqn:prodmetric}
	\tdis{\seqx}{\seqxp}\eqdef\sum_{i=0}^{\infty}\rho^{-i}\dis{\xs{i}}{\xps{i}}.
\end{equation*}
The \deffont{shift map} \selfmap{\sigma}{\XtoN}, defined by 
	\begin{equation*}\label{eqn:shiftmap}
		\sigma(\seqx)_{i}\eqdef\xs{i+1},\quad\inflist{i}{0,1},
	\end{equation*}
is continuous, and the set
	\begin{equation*}\label{eqn:inverselim}
		\Xf\eqdef\setbld{\seqx=\xs{0},\xs{1},...}{\xs{i+1}=\fof{\xs{i}}\text{ for }\inflist{i}{0,1}}
	\end{equation*}
is invariant under the shift; in fact the map $x\mapsto\seqx\eqdef x,\fof{x},\ftoof{2}{x},...$ is a topological
conjugacy between \fXX{} and the restriction of \sigme{} to \Xf{}.  Thus, they have the same topological entropy.  But the entropy of \sigme{} on \Xf{} can be defined by analogy with the topological entropy of a  subshift of the shift on sequences from a finite alphabet---which is given by the growth rate of the number of words of length $n$. Given \epsgo{} and \inNat{n}, define an \emph{\eps-cube} of order $n$ to be the subset of \XtoN{} obtained by specifiying $n$ open \eps-balls $\Bs{0},...,\Bs{n-1}$ and considering all sequences \seqx{} with
$\xs{i}\in\Bs{i}$ for \ilist{i}{0}{n-1} (and no conditions on \xs{i} for $i\geq n$).  Then define the \emph{\eps-capacity} of \Xf{} of order $n$ to be the minimum number of \eps-cubes of order $n$ needed to cover \Xf{};  the entropy of \sigme{} on \Xf{} equals the limit as $\eps\to0$ of the growth rate of the \eps-capacity of \Xf{} of order $n$.

\section{Definitions of Entropy in Non-Compact Spaces}\label{sec:noncompact}

How can the preceding definitions be adapted to \fXX{} when \X{} is \emph{not} compact?

\subsection{Bowen Compacta Entropy}
Bowen's first approach to this question\footnote{In \cite{Bowen:noncompact} he proposed a different answer, based on ideas related to Hausdorff dimension, which we do not consider here.},
in \cite{Bowen:homog}, was to note that the entropy of \f{} on any compact subset $K\subset\X$ is  still well-defined by \refer{eqn}{sepent}.  (One should note, however, that the set $K$ is \emph{not} assumed to be invariant under \f{};  we are thus measuring the dispersion of orbits \emph{emanating from} $K$, but not confined to it, and hence even at this level, the entropy of \f{} ``on'' a compact set $K$ can---and does---depend on the choice of \dist{} we use to calculate it.)  Then he defined the topological entropy of a uniformly continuous \Selfmap{} \fXX{} on an arbitrary metric space to be the supremum of its entropy on all compact subsets of \X{}.  We will distinguish this definition by referring to it as \deffont{Bowen compacta entropy}:
\begin{equation}\label{eqn:compacta}
	\hcompof{\f}{\X}{\dist}\eqdef\sup_{K\subset\X\text{ compact}}\hsepof{\f}{K}{\dist}.
\end{equation}
It is not clear that the assumption of uniform continuity plays any role in this \emph{definition}, although 
the examples motivating Bowen were all uniformly continuous.   In fact, when \f{} is uniformly continuous, then it extends continuously to the completion of \Xd{}, which (assuming the metric space \Xd{} is totally bounded) is compact, and then by \refer{cor}{dense} we are looking at the standard, compact case.
We shall avoid the assumption of uniform continuity for our purposes.  However, \refer{rmk}{bad}, as well as examples outlined by Walters \cite[p. 176]{Walters},
show that when the assumption of uniform continuity is dropped, a number of useful properties  of entropy that hold for uniformly continuous maps can fail.

\subsection{Friedland Entropy}
Rational maps on complex projective space \CP{n} are defined only on a subset of the compact space \CP{n}, and are definitely not uniformly continuous.  To define topological entropy for such maps, Friedland \cite{Friedland:polyrat, Friedland:alg} adapted Gromov's point of view.  Recall that a metric space \Xd{} is \deffont{totally bounded} if for every \epsgo{} there exists a finite cover of \X{} by balls of radius \eps{}.  It is a standard fact that this condition is equivalent to the possibility of embedding \Xd{} isometrically in some compact metric space \cite[p. 198]{Kelley}, \cite[p.276]{Munkres}, which can be taken to be the completion \hXd{} of \Xd{}.  Thus, if \fXX{} is a continuous \Selfmap{} of a totally bounded metric space \Xd{} then we can form the set $\Xf\subset\hX^{\Nat}$ as before;  the space $\hX^{\Nat}$ is compact, and \Xf{} is invariant under the shift map \selfmap{\sigma}{\hX^{\Nat}}, so we can take the topological entropy of \sigme{} on the closure of \Xf{} in \hX{}.  Note that this requires no assumptions of uniform continuity on \f{}.  This gives the \deffont{Friedland entropy}
\begin{equation}\label{eqn:Friedland}
	\hFf\eqdef\htop{\sigma\,|\,\clos\Xf}.
\end{equation}

Friedland's definition does not explicitly involve a metric, but rather an embedding of the space \X{} in some compact topological space.  Of course, if \X{} carries a totally bounded metric, it singles out such an embedding.  It should also be noted that Friedland's formulation was used by him to extend the notion of entropy beyond iterated mappings, to more general relations \cite{Friedland:graphs}. 

\subsection{Bowen-Dinaburg Entropy}
Finally, we consider a straightforward translation of Bowen-Dinaburg entropy to this setting.  Suppose again that \fXX{} is a continuous \Selfmap{} of the metric space \Xd{}.  Note that the condition that the numbers \minspan{\X}{\dist}{\eps} and \ms{\X}{\dist}{\eps}  are finite for all \epsgo{} is precisely total boundedness of \Xd{}, so we assume this.   This implies that the Bowen-Dinaburg metrics \dnf{} are  also totally bounded: the easiest way to see this is to note that the completion of the maximum metric on $\X^{n}$ 
\begin{equation*}
	\diss{max}{\seqx=(\xs{1},...,\xs{n})}{\seqxp=(\xps{1},...,\xps{n})}=\max_{\ilist{i}{1}{n}}\dis{\xs{i}}{\xps{i}}
\end{equation*} 
is the maximum metric on $\hX^{n}$. The latter is compact, and contains the collection of orbit segments $\seqx=(x,\fofx,...,\ftoof{n-1}{x})$.  In particular, the numbers \minspan{\X}{\dnf}{\eps} and \ms{\X}{\dnf}{\eps} are finite for \inflist{n}{0,1,2} and all \epsgo{}.  This means that the definition of \hsepof{\f}{K}{\dist} given by Equations \eqref{eqn:BowenDina} and \eqref{eqn:sepent} makes sense for any continuous (not necessarily uniformly continuous) \Selfmap{} \fXX{} and any subset $K\subset\X$ whenever the metric \dist{} on the space \X{} is totally bounded.  We will be interested primarily in the case $K=\X$.

\subsection{Coding}

In some contexts, there is a natural way of replacing a noncompact
dynamical system by a symbolic system, that is to say a subshift
of the full shift over some finite alphabet.  Interval exchange
transformations provide a simple example.  Consider for instance 
the interval exchange on the open interval $(0,1)$ that maps
$(0,\alpha)$ affinely to $(1-\alpha,1)$ and $(\alpha,1)$ 
affinely to $(0,1-\alpha)$, with $\alpha$ irrational.  (Maps like 
this arise as Poincar\'e sections for the billiards flow on a 
compact table with internal corners; the endpoints of the intervals 
being exchanged correspond to orbit-segments that hit a corner.)  
This map is not defined at points $x$ in $(0,1)$ that are fractional 
parts of some multiple of $\alpha$; therefore, if we want to view 
the dynamics as arising from iteration of a function from a
domain to itself, that domain should be $(0,1)$ minus a countable
dense set of points.  Call this non-compact domain $X$.  We can 
map our dynamical system on $X$ into the 2-shift by partitioning 
$X$ into $X \cap (0,\alpha)$ and $X \cap (\alpha,1)$ and then
coding an orbit in $X$ by a binary string in the usual way.  Let
$\psi$ denote the map from $X$ to the 2-shift.  When we take the 
closure of $\psi(X)$, we add countably many limit points, obtaining 
a compact shift-invariant set $X'$.  If one can show that the ergodic 
non-atomic shift-invariant measures on $X'$ are all supported on
$\psi(X)$ and hence correspond to the ergodic non-atomic invariant 
measures on $X$ (as is the case for systems derived from polygonal 
billiards \cite{Kat87}), then one may feel justified in regarding 
the two systems as closely linked and defining the entropy of the 
former to equal the entropy of the latter.  We mention this approach 
to compactifying dynamical systems,
but we will not pursue it beyond suggesting that the affinity 
with Friedland's approach ought to be explored further.

\section{Relation between Entropies}\label{sec:relation}

We need to address two issues with respect to the three definitions of entropy in \refer{sec}{noncompact}.
First, we need to establish the extent to which they are invariants, and second, we need to establish the relations between them.

As noted in \refer{sec}{compact}, when \X{} is compact, all three definitions yield the same quantity, which agrees with the Adler-Konheim-McAndrew entropy, and the latter is clearly an invariant of topological conjugacy. (This is because a homeomorphism preserves open covers and their cardinality.)  However, in the context of totally bounded spaces, \hAKMf{} does not in general agree with any of these quantities (we will see this in \refer{sec}{escape}), so we need to attack the invariance question differently.

Suppose \fXd{} and \tfXd{} are (not necessarily uniformly) continuous \Selfmap{}s of totally bounded metric spaces.  A \deffont{semiconjugacy}\footnote{Sometimes \f{} is referred to as a \deffont{factor} of \tilf{} in this case.} from \tfXd{} to \fXd{} is a continuous surjection 
\map{h}{\pairord{\tX}{\tdist}}{\pairord{\X}{\dist}} satisfying $\comp{h}{\tilf}=\comp{\f}{h}$; it is a \deffont{uniform semiconjugacy} if $h$ is uniformly continuous with respect to the metrics \tdist{} and \dist{}: 
that is, for each \epsgo{}, there exists 
$\teps>0$ such that  $\tdis{x}{\xp}<\teps$ implies $\dis{\hof{x}}{\hof{\xp}}<\eps$.  The choice of \teps{} given \eps{} is a \deffont{modulus of continuity} for $h$.
When $h$ is a uniformly continuous homeomorphism with uniformly continuous inverse, it is a \deffont{uniform conjugacy}.

\begin{lemma}\label{lem:semiconj} Suppose \tfXd{} and \fXd{} are continuous \Selfmap{}s of
	totally bounded metric spaces.
	\begin{enumerate}
		\item A uniform semiconjugacy $h$ from \tfXd{} to \fXd{} is uniformly continuous 
			with respect to each pair of corresponding Bowen-Dinaburg metrics, 
			with the same modulus of continuity.
		\item In particular, for each \inNat{n} and \epsgo{},
			\begin{align}
				\ms{\tX}{\tdBowen{\tilf}{n}}{\teps}&\geq\ms{\X}{\dBowen{f}{n}}{\eps}
					\label{eqn:semiconjms}\\
				\minsp{\tX}{\tdBowen{\tilf}{n}}{\teps}&\geq\minsp{\X}{\dBowen{f}{n}}{\eps}.
					\label{eqn:semiconjminsp}
			\end{align}
		\item It follows that
			\begin{equation}\label{eqn:semiconjhsep}
				\hsepof{\tilf}{\tX}{\tdist}\geq\hsepof{\f}{\X}{\dist}.
			\end{equation}
		\item Thus \hsep{} is an invariant of uniform conjugacy.
	\end{enumerate}
\end{lemma}

\begin{proof}
	To establish the first statement, let \eps{} and \teps{} be related 
	as in the definition of uniform continuity for $h$, as above.  
	Note that if $\tdB{\tilf}{n}{x}{\xp}<\teps$ then for each \ilist{i}{0}{n-1}
	we have $\tdis{\toof{\tilf}{i}{x}}{\toof{\tilf}{i}{\xp}}<\teps$, which implies (for each $i$) that 
	$\dis{\hof{\toof{\tilf}{i}{x}}}{\hof{\toof{\tilf}{i}{\xp}}}<\eps$, 
	but since $\comp{h}{\too{\tilf}{i}}=\comp{\fto{i}}{h}$, this is the same as 
	$\dnfof{\hof{x}}{\hof{\xp}}<\eps$.
	
	If $S\subset\X$ is \eps-separated with respect to \dBowen{\f}{n}, form $\tilde{S}\subset\tX$ 
	by picking a single preimage of each element of $S$.  
	Then $\tilde{S}$ must be \teps-separated with respect to \tdBowen{\tilf}{n}, by part (1),
	and has the same cardinality as $S$,
	giving \refer{eqn}{semiconjms}.
	
	If $\tilde{S}\subset\tX$ is \teps-spanning with respect to \tdBowen{\tilf}{n},  
	then $S=\hof{\tilde{S}}\subset\X$
	is \eps-spanning with respect to \dBowen{\f}{n}, since (by surjectivity of $h$) 
	given $x\in\X$ we can pick $\tilde{x}\in\tX$ with 
	$\hof{\tilde{x}}=x$, and then pick $\tilde{s}\in\tilde{S}$ with $\tdB{\tilf}{n}{\tilde{x}}{\tilde{s}}<\teps$, 
	and it follows
	that $s=\hof{\tilde{s}}\in S$ satisfies $\dB{\f}{n}{x}{s}<\eps$.  \refer{eqn}{semiconjminsp} follows.
	
	The other two statements are immediate consequences.
\end{proof}

It is an immediate consequence that Bowen compacta entropy is an invariant of uniform conjugacy;  however, to obtain the corresponding analogue of \refer{eqn}{semiconjhsep} we need to know that the preimage of every compact subset $K\subset\X$ is a compact subset of \tX{}---that is, the semiconjugacy must be \deffont{proper}.  With this assumption, we have

\begin{remark}\label{rmk:semiconj}
	If \map{h}{\pairord{\tX}{\tdist}}{\pairord{\X}{\dist}} is a proper, uniform semiconjugacy from \tfXX{}
	to \fXX{}, then
		\begin{equation*}
			\hcompof{\tilf}{\tX}{\tdist}\geq\hcompof{\f}{\X}{\dist}.
		\end{equation*}
\end{remark}

While we could establish directly that Friedland entropy is an invariant of uniform conjugacy, this will be an immediate corollary of the equality between it and Bowen-Dinaburg entropy, which we establish next, based on two observations.   The first observation is

\begin{lemma}\label{lem:densesep}
	Suppose $\Y\subset\X$ is a dense subset of the totally bounded metric space \Xd{}.
	
	Then
	\begin{enumerate}
		\item For any $\peps$, $0<\peps<\eps$,
			$$\ms{\Y}{\dist}{\peps}\geq\ms{\X}{\dist}{\eps}.$$
		\item For any $\peps>\eps$,
			$$\minsp{\Y}{\dist}{\peps}\leq\minsp{\X}{\dist}{\eps}.$$
	\end{enumerate}
\end{lemma}

\begin{proof}
	To see the first statement, suppose $S\subset\X$ is \eps-separated.  Given $0<\peps<\eps$,
	let $\delta=\half(\eps-\peps)>0$. For each $s\in S$, pick $y\in\Y$ with $\dis{s}{y}<\delta$, and
	let $S\pr\subset\Y$ be the resulting set of $y$'s.  For $s\neq \spr\in S$, let $y,\yp$ be the
	corresponding points in $S\pr$.  Then the triangle inequality gives
		$$\dis{y}{\yp}\geq\dis{s}{\spr}-[\dis{s}{y}+\dis{\spr}{\yp}]>\eps-2\delta=\peps$$
	so $S\pr\subset\Y$ is \peps-separated and has cardinality the same as $S$.
	
	To see the second statement, suppose $S\subset\X$ \eps-spans \X{}.  Given $\peps>\eps$,
	let $\delta=\peps-\eps$, and for each $s\in S$ pick $\spr\in\Y$ with $\dis{\spr}{s}<\delta$.
	Again, form the set $S\pr\subset\Y$ of all such \spr{}'s.  If $y\in\Y\subset\X$, we can find $s\in S$ with
	$\dis{s}{y}<\eps$;  then 
		$$\dis{\spr}{y}\leq\dis{\spr}{s}+\dis{s}{y}<\delta+\eps=\peps.$$
	Thus, $S\pr\subset\Y$ is a \peps-spanning subset of \Y{} with cardinality at most that of $S$.
	
\end{proof}

\begin{corollary}\label{cor:dense}
	Suppose \fXd{} is a continuous self-map of a totally bounded metric space and $\Y\subset\X$
	is a dense subset.  Then
		$$\hsepof{\f}{\X}{\dist}=\hsepof{\f}{\Y}{\dist}.$$
\end{corollary}

Now, consider the map \map{\vphi}{\X}{\Xf} taking \inside{x}{\X} to its orbit
\begin{equation*}
	\vphiof{x}=\seqx\eqdef x,\fof{x},\ftoof{1}{x},....
\end{equation*}
This is clearly a bijection, whose inverse is the restriction to \Xf{} of the projection \map{\pi}{\XtoN}{\X} to the first factor of \XtoN{}:
\begin{equation*}
	\piof{\seqx\eqdef\xs{0},\xs{1},...}=\xs{0}.
\end{equation*}
Furthermore, these maps are equivariant with respect to \f{} and the shift map \sigme{}, and in particular \pie{} is a semiconjugacy from \selfmap{\sigma}{\Xf} to \fXX{}.

\begin{lemma}\label{lem:orbit}

	\begin{enumerate}
		\item For each $n\in\Nat$ the projection \pie{} is a uniformly continuous 
		semiconjugacy from \hXfd{} to \Xdnf{}.
		\item\label{lem:orbit2} For each \epsgo{}, there exists \inNat{N} such that for all $n\geq N$,
			$$\dB{f}{n}{x}{\xp}<\eps\imply\hdis{\vphiof{x}}{\vphiof{\xp}}<(N+1)\eps.$$
	\end{enumerate}
	
\end{lemma}

\begin{proof}
\emph{(1)}  The projection $\seqx\mapsto\xs{0}$ is uniformly continuous from \hXfd{} to \Xdnf{}, 
	with modulus of continuity $\peps=\rho^{-n}\eps$: if $\hdis{x}{\xp}<\peps$, then for each $j<n$
	\begin{equation*}
		\rho^{-j}\dis{\ftoof{j}{x}}{\ftoof{j}{\xp}}
			\leq\sum_{i=0}^{\infty}\rho^{-i}\dis{\ftoof{i}{x}}{\ftoof{i}{\xp}}
			=\hdis{x}{\xp}<\peps
	\end{equation*}
	or
	\begin{equation*}
		\dis{\ftoof{j}{x}}{\ftoof{j}{\xp}}<\rho^{j}\peps<\eps,
	\end{equation*}
	because $j<n$ and $\rho>1$. Maximizing over $j<n$ gives 
	\begin{equation*}
		\dnfof{x}{\xp}<\eps
	\end{equation*}
	as required.
		
	\emph{(2)} Suppose we are given \epsgo{}.  Since \X{} is totally bounded, it has finite
	diameter, and we can
	find \inNat{N} so that
		$$\sum_{i=N}^{\infty}\rho^{-i}\diam\Xd<\eps.$$
	Now, for $n\geq N$, if $\dB{f}{n}{x}{\xp}<\eps$, then
		\begin{align*}
			\hdis{\vphiof{x}}{\vphiof{\xp}}
			&=\sum_{i=0}^{N-1}\rho^{-i}\dis{\ftoof{i}{x}}{\ftoof{i}{\xp}}
			+\sum_{i=N}^{\infty}\rho^{-i}\dis{\ftoof{i}{x}}{\ftoof{i}{\xp}}\\
			&\leq N\max_{i<N}\dis{\ftoof{i}{x}}{\ftoof{i}{\xp}}
			+\sum_{i=N}^{\infty}\rho^{-i}\diam\Xd\\
			&\leq N\dB{f}{n}{x}{\xp}+\eps\\
			&<(N+1)\eps
		\end{align*}
	as required.
\end{proof} 

Using these observations, we can prove

\begin{theorem}
	For every continuous map \fXd{} of a totally bounded metric space to itself,
		$$\hFf=\hsepof{\f}{\X}{\dist}\geq\hcompof{\f}{\X}{\dist}.$$
\end{theorem}

\begin{proof}
	First, we establish the equality between Friedland entropy and Bowen-Dinaburg entropy.
	
	By \refer{eqn}{semiconjhsep} in \refer{lem}{semiconj}, the first statement of \refer{lem}{orbit}
	shows that 
		$$\hsepof{\f}{\X}{\dist}\leq\hsepof{\sigma}{\Xf}{\hdist}$$
	while \refer{cor}{dense} shows that 
		$$\hsepof{\sigma}{\Xf}{\hdist}=\hsepof{\sigma}{\clos\Xf}{\cdist}=:\hFf.$$
		
	In the other direction, given \epsgo{} pick $N$ as in \refer{lem}{orbit}\eqref{lem:orbit2}.
	Then for every $n>N$, we see that for every set $S\subset\X$ which is
	$\frac{\eps}{N+1}$-spanning with respect to \dnf{},
	\vphiof{S} \eps-spans \Xf{} with respect to \hdist{}.  It follows that (again for $n>N$)
		$$\minsp{\X}{\dnf}{\frac{\eps}{N+1}}\geq\minsp{\Xf}{\hdist}{\eps}.$$
	But then, fixing \eps{}, we have
		$$\GR{\minsp{\X}{\dnf}{\frac{\eps}{N+1}}}\geq\GR{\minsp{\Xf}{\hdist}{\eps}}$$
	and taking the limit as $\eps\to0$ we get
		$$\hsepof{\f}{\X}{\dist}\geq\hFf.$$
		
	Finally, to show that Bowen-Dinaburg entropy dominates Bowen compacta entropy, it is enough to note that
	the latter is the supremum of Bowen-Dinaburg entropy over the collection of compact subsets of \X{},
	and that the Bowen-Dinaburg entropy of the restriction of a map to a subset 
	is at most that of the map on the ambient space.
\end{proof}

\section{Dependence on uniform structure}\label{sec:escape}

The original definition of Adler-Konheim-McAndrew made no reference to a metric structure, but was formulated in the context  of compact topological spaces.  It might be tempting, therefore, to adopt \refer{eqn}{AKM} as a definition of topological entropy (that is, using only \emph{finite} open covers, even if the space is not compact).  While this conjugacy invariant is independent of any metric or uniform structure, it is unfortunately also infinite for most essentially non-compact examples.

In a non-compact space, an orbit can \emph{escape to infinity} in the sense that it eventually leaves every compact subset, or equivalently, it has no accumulation points (i.e, $\omega(x)=\emptyset$).

\begin{prop}\label{prop:escape} Suppose \fXX{}, a continuous \Selfmap{} of a topological space, 
has an orbit that escapes to infinity.  Then:
\begin{enumerate}
	\item $\hAKMf{}=\infty$.
	\item If \X{} carries a totally bounded metric \X{}, then for each \inNat{N} 
	there is a totally bounded metric \dn{} on \X{}
	such that
	\begin{equation*}
		\hsepof{f}{\X}{\dn}\geq N.
	\end{equation*}
\end{enumerate}

\end{prop}

\begin{proof}
	First, we calculate the Adler-Konheim-McAndrew entropy of $f$.  Suppose the orbit \infseq{\ftoof{i}{x}}{i}{0} of $x$ has no accumulation points;  then as a set it is closed, countable and discrete.  Pick \Us{i}, \inflist{i}{0,1} disjoint open neighborhoods of \ftoof{i}{x}.  Note that $V\eqdef\X\setminus\setbld{\ftoof{i}{x}}{\inflist{i}{0,1}}$ is open.
	
	Given a positive integer $N$, let \infseq{s_{i}}{i}{0} be a sequence of integers 
	between $0$ and $N-1$ ($s_{i}\in\single{0,...,N-1}$
	for all $i$) containing every finite ``word'' $w=\ws{0}...\ws{n-1}$, in the
	sense that for some \is{w}, $s_{i_{w}+j}=\ws{j}$, \ilist{i}{0}{N-1}.  
	Now, group the neighborhoods \Us{i} of \ftoof{i}{x} according to the values of $s_{i}$:  
	for \ilist{k}{1}{N-1}, define
	\begin{equation*}
		\As{k}\eqdef\bigcup\setbld{\Us{i}}{s_{i}=k},\quad\ilist{k}{1}{N-1}
	\end{equation*}
	and define \As{0} similarly, but adjoining $V$
	\begin{equation*}
		\As{0}\eqdef V\cup\bigcup\setbld{\Us{i}}{s_{i}=0}.
	\end{equation*}
	
	Then 
	\begin{equation*}
		\alpha\eqdef\single{\As{0},...,\As{N-1}}
	\end{equation*}
	is an open cover of \X{} by $N$ open sets.  It has no proper subcover, because each \ftoof{i}{x}
	belongs only to \As{s_{i}}.  Furthermore, if for each ``word'' $w=\ws{0}...\ws{n-1}$ 
	we define an element of \alfs{n} by
	\begin{equation*}
		\As{w}\eqdef\As{\ws{0}}\cap\fpreof{1}{\As{\ws{1}}}\cap\dots\cap\fpreof{(n+1)}{\As{\ws{n-1}}}
	\end{equation*}
	then the point \ftoof{\is{w}}{x} belongs to \As{w} and to none of the sets 
	corresponding to other words of length $n$.  This shows that the refinement \alfs{n} 
	also has no proper subcovers;  thus $\Hof{\alfs{n}}=n\log N$ so
	\begin{equation*}
		\hof{\f,\alpha}=\log N
	\end{equation*}
	and the supremum over all finite covers is
	\begin{equation*}
		\hAKMf{}=\infty,
	\end{equation*}
	establishing (1).
	
	Now suppose \X{} is embeddable in a compact metric space, or equivalently, that we can find
	a totally bounded metric $d$ on \X{}.  Fix $N$, \Us{i}, and $s_{i}$ as above.  Using the Urysohn
	Lemma, let \phie{} be a continuous bounded function supported inside the neighborhoods \Us{i}
	of the orbit such that
	\begin{equation*}
		\phiof{\ftoof{i}{x}}=s_{i}+1,\quad \inflist{i}{0,1}.
	\end{equation*}
	
	Then \X{} embeds in $X\times\Reals$ as the graph of \phie{}, and we can use the natural product
	metric, $\tilde{d}$, which is totally bounded on the graph.  It is easy to see that if $s_{i}\neq s_{j}$
	then $\tilde{d}(\ftoof{s_{i}}{x},\ftoof{s_{j}}{x})\geq 1$;  from this it follows that by picking a point
	of the orbit in each of the sets $\As{w}\in\alfs{n}$ we obtain an $(n,1)$-separated set of cardinality
	$N^{n}$, and the rest follows by standard arguments.
	
\end{proof}

It is easy to see that a similar phenomenon occurs if some point \xs{0} ``escapes to infinity'' in backward time, in the sense that there exists a sequence of successive preimages \xs{-i}, \inflist{i}{0,1}, ($\fof{\xs{-i}}=\xs{1-i}$ for $i\geq1$) with no accumulation points.

\refer{prop}{escape} shows that we cannot hope to avoid the effects of a choice of metric on the definition of topological entropy in a non-compact setting.  An elaboration of the technique of proof for this result can be used to illustrate the ways that our two metric notions of entropy---as well as the relation between them---can be affected by the choice of metric. 

In general, one way to define a totally bounded metric on a topological space \X{} 
is to (topologically) embed \X{} in a compact metric space; then the restriction of the metric to the 
embedded image of \X{} defines a metric on \X{} which is clearly totally bounded.  

Using this trick, we
can establish
\begin{prop}\label{prop:arbent}
	Suppose \selfmap{\f}{\opint{0}{1}} is a homeomorphism satisfying 
	\begin{equation*}
		\fofx>x\text{ for all }x\in\opint{0}{1}.
	\end{equation*} 
	
	Then for any
	positive integer $N$ there exists a totally bounded metric \dn{} on \opint{0}{1} such that
	\begin{align}
		\hsepof{\f}{\opint{0}{1}}{\dn}=\hsepof{\finv}{\opint{0}{1}}{\dn}&=\log\, N
			\label{eqn:arbcomponent0}\\
		\hcompof{\f}{\opint{0}{1}}{\dn}&=\log\,N\label{eqn:arbcompent1}\\
		\intertext{and}
		\hcompof{\finv}{\opint{0}{1}}{\dn}&=0\label{eqn:arbcompent2}.
	\end{align}
\end{prop}

Note that by \refer{prop}{escape}, we already know that 
\begin{equation*}
	\hAKMf{}=\hAKMof{\finv}=\infty.
\end{equation*}

We note in passing that this result is more general than may at first appear, because of the
following observation, which is a standard exercise in basic dynamics.

\begin{remark}\label{rmk:fixfree}
	Any two fixedpoint-free self-homeomorphisms of open intervals are topologically conjugate.
\end{remark}
 
 In particular, \f{} and \finv{} are topologically conjugate, and hence the metric for which 
 $\hcompof{\f}{\opint{0}{1}}{\dist}=\log N$ and $\hcompof{\finv}{\opint{0}{1}}{\dn}=0$ and can be turned into a metric \dist{} for which 
 $\hcompof{\f}{\opint{0}{1}}{\dist}=0$  and $\hcompof{\finv}{\opint{0}{1}}{\dist}=\log N$ by composing the embedding defining \dn{} with the homeomorphism that conjugates \f{} with \finv{}.
 
\begin{proof}[Proof of \refer{prop}{arbent}]
	In the standard metric on \opint{0}{1}, \f{} automatically extends to the closed interval by fixing the
	endpoints.  This extension has topological (hence Bowen-Dinaburg) entropy zero, 
	since the topological entropy of any continuous \Selfmap{} of a compact space equals the entropy
	of its restriction to its nonwandering set (\cite{Bowen:topaxA, Xiong:nonwandering}), 
	which in this case consists of the two endpoints.  But then
	\refer{cor}{dense} tells us that the same holds for \selfmap{\f}{\opint{0}{1}}.  This takes care
	of the special case $N=1$ of the lemma.
	
	For $N\geq2$, we will construct a family of examples depending on $N$.  
	These all come from embedding \opint{0}{1} in a piecewise-linear
	version of the ``topologist's sine curve'' and restricting the standard metric on \Realsto{2} 
	to this curve.  Given a strictly decreasing sequence  \infseq{\as{k}}{k}{0} 
	with $\as{0}=1$ 
	and $\lim\as{k}=0$, define the piecewise-affine 
	function \Phie{} by
	\begin{equation*}
		\Phiof{\as{k}}=(-1)^{k}
	\end{equation*}
	making \Phie{} affine on each interval $\Is{k}=\clint{\as{k+1}}{\as{k}}$.  
	We denote the graph of \Phie{} on \opint{0}{1} by
	\begin{equation*}
		\Gam\eqdef\setbld{\pairpar{x}{\Phiof{x}}}{0<x<1}.
	\end{equation*}
	Note that the closure of \Gam{} is the compact subset of \Realsto{2}
	consisting of \Gam{},  
	the point \pairpar{1}{1},
	and the vertical interval \setbld{\pairpar{0}{t}}{-1\leq t\leq 1}.  We will choose the sequence 
	\infseq{\as{k}}{k}{0} at our convenience during the construction.
	
	We will construct an example satisfying 
	\refer{eqn}{arbcompent1} whose inverse will satisfy 
	\refer{eqn}{arbcompent2}.  (Note that by \refer{rmk}{fixfree} \f{} and
	\inverse{f} are topologically conjugate, so if the metric coming from one embedding of 
	\opint{0}{1} into a compact space gives Bowen compacta entropy zero for \inverse{f}, 
	then composing
	this embedding with a conjugacy between \f{} and \inverse{f} will yield a different embedding
	for which the corresponding metric will give \f{} itself zero Bowen compacta entropy.)
	
	First we construct a model of
	\f{} on \opint{0}{1} by defining \selfmap{\g}{\clint{0}{1}} to fix $0$ and $1$, and make the sequence
	\infseq{\frac{1}{j}}{j}{2} a (forward) orbit
	\begin{equation*}
		\gof{\frac{1}{j}}=\frac{1}{j+1},\quad \inflist{j}{2,3}
	\end{equation*}
	and requiring \g{} to be affine on each of the intervals \clint{\frac{1}{j+1}}{\frac{1}{j}}, \inflist{j}{0,1}.
	By \refer{rmk}{fixfree}, $\g|\opint{0}{1}$ is conjugate to our original map \f{}.
	
	Next, we construct the sequence \single{\as{k}} by setting $\as{0}=1$, $\as{1}=\half$, 
	and then dividing each
	of the intervals \clint{\recip{j+1}}{\recip{j}}, \inflist{j}{2,3} into $N^{j-1}$ intervals of equal length 
	by successive points of the sequence \single{\as{k}}; this means that 
	\begin{equation*}
		\recip{\ell+2}=\as{1+N+...+N^{\ell}}\quad\text{ for }\inflist{\ell}{0,1}.
	\end{equation*}
	The affine property of \g{} implies that each interval $\Is{k}=\clint{\as{k+1}}{\as{k}}$, 
	\inflist{k}{1} maps
	to a union of $N$ such intervals.
	
	Now, if we use the homeomorphism between \opint{0}{1} and \Gam{} 
	defined by $x\mapsto(x,\Phiof{x})$,
	we conjugate \g{} with a homeomorphism \selfmap{G}{\Gam} that takes each of the ``laps'' 
	of \Phie{} 
	(the pieces $\Gams{k}\eqdef\setbld{(x,\Phiof{x})}{x\in\Is{k}}$ of \Gam{} joining $y=-1$ to $y=1$) 
	other than the rightmost one and ``crumples'' it
	into $N$ laps.  
	
	\begin{claim}\label{claim:good}
		For \dn{} the metric inherited from the embedding of \Gam{} in \Realsto{2},
		\begin{equation}\label{eqn:BDBc}
			\hsepof{G}{\opint{0}{1}}{\dn}=\hcompof{G}{\opint{0}{1}}{\dn}=\log N.
		\end{equation}
	\end{claim}

	To see this, consider the iterates of $G$
	on some lap \Gams{k} of \Gam{}.  The first iterate maps
	\Gams{k} onto a union of $N$ sets of the form \Gams{k\pr+i}, \ilist{i}{0}{N-1}, 
	and in particular for $(x,y)\in\Gams{k}$,  $G(x,y)=(\gof{x},\hof{y})$, where \selfmap{h}{\clint{-1}{1}}
	maps each of the $N$ subintervals \clint{-1+\frac{2j}{N}}{-1+\frac{2(j+1)}{N}} (\ilist{j}{0}{N-1}) 
	affinely onto
	\clint{-1}{1}.  It is well-known (\eg{} \cite{Misiurewicz-Szlenk}, \cite[p. 205]{ALM}) 
	that $h$ has entropy $\log N$,
	and hence has \nepssep{} sets of cardinality growing like $N^{n}$ for any \epsgo{}.  If we pick
	such a set of $y$-values, then the corresponding points of \Gams{k} form an \nepssep{} set for $G$.
	This shows that both the Bowen-Dinaburg and Bowen compacta entropies of $G$ 
	are at least $\log N$;  but the 
	reverse inequality can be seen from the fact that the ``base'' map \g{} has entropy zero.

	\begin{claim}\label{claim:bad}
		For \dn{} the metric inherited from the embedding of \Gam{} in \Realsto{2},
		\begin{equation}\label{eqn:BDinv}
			\hsepof{\inverse{G}}{\opint{0}{1}}{\dn}=\log N
		\end{equation}
		but
		\begin{equation}\label{eqn:Bcinv}
			\hcompof{\inverse{G}}{\opint{0}{1}}{\dn}=0.
		\end{equation}
	\end{claim}

	\refer{eqn}{BDinv}
	follows from the standard observation that 
	for any homeomorphism $h$, if $S$ is \nepssep{} under $h$ 
	then $h^{n}(S)$ is \nepssep{} under \inverse{h}.  
	
	However, this argument fails for Bowen compacta entropy,
	because for any compact set $K\subset\Gam$, $G^{n}(K)$ is disjoint from $K$ 
	for $n$ sufficiently large.  
	
	In fact, for large $n$, all images of $K$ under $G^{n}$ will be contained
	in the part of \Gam{} over an arbitrarily small interval with right endpoint $1$, and such a set has
	very small diameter. Thus to be \nepssep{}, a subset of $K$ must already be separated under a
	bounded number of iterates, and the cardinality of such a subset 
	is uniformly bounded independent
	of $n$, and in particular has growth rate zero, establishing \refer{eqn}{Bcinv}.  
	
	Another way to see \refer{eqn}{Bcinv} is to note that $G$ extends continuously to the 
	``right endpoint'' \pairpar{1}{1} of \Gam{}, and the part of \Gam{} lying above any interval
	\clint{a}{1}, with $a>0$, is mapped into itself by this extension. Since its nonwandering set
	is just the right endpoint, the topological entropy of the restriction to this set is zero.  Since any 
	compact subset of \Gam{} is contained in such a piece, the Bowen-Dinaburg entropy of $G$
	on any compact subset is zero, and so the Bowen compacta entropy of $G$ is also zero.
\end{proof}

\begin{remark}\label{rmk:bad}
	The argument for \refer{eqn}{BDinv} shows that the Friedland-Bowen-Dinaburg entropy 
	of a homeomorphism equals that of its inverse (as is the case for topological entropy in compact
	spaces), 
	while Equations \eqref{eqn:BDBc} and \eqref{eqn:Bcinv} show that this fails for Bowen compacta
	entropy.
\end{remark}

Several other examples, showing that Bowen compacta entropy for a non-uniformly continuous map fails to enjoy several useful properties of topological entropy on compact spaces (including that noted above) 
are outlined by Walters \cite[p. 176]{Walters}.

Using \refer{prop}{arbent}, we can construct another family of examples.  

Let 
	$$\X=\opint{0}{1}\times S^{1}$$
be the product of an open interval with the circle $S^{1}\eqdef\Reals/\Integers$.  Topologically, \X{} is an open annulus.  Define \fXX{} by
	\begin{equation*}
		\fof{x,\theta}=(x^{2},2\theta\text{ mod }\Integers).
	\end{equation*}

The embedding \embedding{\varphi}{\X}{\disc}, of \X{} into the open unit disc 
$\disc\eqdef\setbld{z\in\Cx}{\abs{z}\leq1}$
given by 
	$$\varphi(x,\theta)=x e^{2\pi i\theta}$$
conjugates \f{} with the restriction of the quadratic map $q:z\mapsto z^{2}$ to the punctured open unit disc $\pdisc{}=\setbld{z\in\Cx}{0<\abs{z}<1}$.  From the Friedland version, it is clear that, with respect to the (pullback to \X{} via \vphi{} of the) standard metric 
\dist{} on \disc{}, 
	$$\hsepof{\f}{\X}{\dist}=\htop{q}=\log 2.$$
But any compact set $K\subset\X$ embeds in some closed disc $B_{r}(0)$ of radius $r<1$, on which the topological entropy of $q$ is zero (the nonwandering set is just the fixedpoint at $z=0$), so that for such a set
	$$\hsepof{\f}{K}{\dist}=\htop{q|B_{r}(0)}=0.$$
	
Now, if instead we embed \X{} in \disc{} via
	$$(x,\theta)\mapsto \left((1-x)e^{2\pi i\theta}\right)$$
the only difference is that under \f{} the radial component increases toward $1$ instead of decreasing toward $0$.
The Bowen-Dinaburg entropy remains at $\log 2$, but this time for $K$ any circle centered at $0$ in \disc{} we
can find some \inNat{N} such that for all $n>N$ \ftoof{n}{K} is a circle very close to the boundary circle 
\single{\abs{z}=1}, and so the Bowen compacta entropy is also $\log 2$.

If we embed \X{} in the unit sphere $S^{2}\subset\Realsto{3}$ with its standard metric, by identifying each ``boundary circle'' of \X{} to a point, we find that both the Bowen-Dinaburg entropy and the Bowen compacta entropy are zero, since the induced map on the sphere has as its nonwandering set the two fixedpoints identified with the boundary circles of \X{}.

Finally, consider the metric \dist{} on \pdisc{} obtained by applying the construction of \refer{prop}{arbent} to the radial component (the disc ``crumples'' infinitely often near the origin, which itself is not part of the disc).  Then the arguments in the proof of \refer{prop}{arbent} show that we can achieve
\begin{equation*}
	\hsepof{f}{\X}{\dist}=\hcompof{f}{\X}{\dist}=\log N
\end{equation*}
for any integer $N\geq2$, and, if we ``crumple'' near the outer boundary instead of near the origin, then the Bowen-Dinaburg entropy will be $\log N$ while the Bowen compacta entropy will be zero.

\bibliography{entropy}
\bibliographystyle{amsalpha}

\end{document}